\documentclass[12pt]{article}
\usepackage{amsmath}
\usepackage{amsthm}
\usepackage{amssymb}
\usepackage{CJK}
\usepackage{vatola}

\def\q{\quad}
\def\e{\equiv}

\def\t{\hbox}
\def\mod#1{\ (\hbox{\rm mod}\ #1)}

\def\a{\alpha}
\def\b{\binom}
\def\s{\sum}
\def\i{\infty}
\def\f{\frac}
\def\te{\text{e}}
\let\pro=\proclaim
\let\endpro=\endproclaim
\def\qtq#1{\q\t{#1}\q}
\theoremstyle{plain}

\TagsOnRight \theoremstyle{definition}

\theoremstyle{lemma}

\begin{document}
\begin{CJK*}{GBK}{kai}

\par\q\par
 \centerline{\bf\large Congruences for sequences analogous to Euler
 numbers} \par\q\par \centerline{Zhi-Hong Sun}
\centerline{School of Mathematical Sciences, Huaiyin Normal
University,} \centerline{Huaian 223001, PR China}
\centerline{E-mail: zhihongsun@yahoo.com} \centerline{Homepage:
http://www.hytc.edu.cn/xsjl/szh} \par\q \newline \centerline{Hai-Yan
Wang}
 \centerline{School of Mathematics and Statistics, Jiangsu Normal
 University,}
 \centerline{Xuzhou 221116, PR China}
 \centerline{E-mail:  2438219717@qq.com}

 \abstract{For a given real number $a$ we define the sequence
$\{E_{n,a}\}$ by
 $E_{0,a}=1$ and $E_{n,a}=-a\sum_{k=1}^{[n/2]}
 \binom n{2k}E_{n-2k,a}$ $(n\ge 1)$,
 where $[x]$ is the greatest integer not exceeding $x$.
  Since $E_{n,1}=E_n$ is the n-th Euler number, $E_{n,a}$ can
   be viewed as a natural
generalization of Euler numbers. In this paper we deduce some
identities and an inversion formula involving $\{E_{n,a}\}$, and
establish congruences for $E_{2n,a}\mod{2^{{\rm ord}_2n+8}}$,
$E_{2n,a}\pmod{3^{{\rm ord}_3n+5}}$ and $E_{2n,a}\pmod{5^{{\rm
ord}_5n+4}}$ provided that $a$ is a nonzero integer, where ${\rm
ord}_pn$ is the least nonnegative integer $\alpha$ such that
$p^{\a}\mid n$ but $p^{\a+1}\nmid n$.
\par\q
\newline MSC: 11B68, 11A07
\newline Keywords: Congruence, Euler number
}
 \endabstract
\section*{1. Introduction }
\par\q The famous Euler numbers  $\{E_n\}$ are  given by
$$E_0=1\q\t{and}\q E_n=-\sum_{k=1}^{[n/2]}\binom n{2k}E_{n-2k}(n\geq1),$$
where $[x]$ is the greatest integer not exceeding $x$. Euler numbers
have many properties and applications. See for example [1-6]. In
[7,8] the first author introduced and studied the similar sequence
$\{U_n\}$  given by
$$U_0=1\q\t{and}\q U_n=-2\sum_{k=1}^{[n/2]}\binom n{2k}U_{n-2k}(n\geq1).$$

\par In this paper, for a given real number $a$ we define the sequence
$\{E_{n,a}\}$ by
 $$E_{0,a}=1\qtq{and}
E_{n,a}=-a\sum_{k=1}^{[n/2]}\binom n{2k}E_{n-2k,a}\ (n\ge 1).\tag
1.1$$ Since $E_{n,1}=E_n$ and $E_{n,2}=U_n$, $E_{n,a}$ can be viewed
as a natural generalization of Euler numbers. Thus, it is
interesting to investigate the properties of $\{E_{n,a}\}$.
\par
 The first few $E_{n,a}$ are as follows:
$$\aligned &E_{2,a}=-a,\ E_{4,a}=-a+6a^2,\
 E_{6,a}=-a+30a^2-90a^3,\\& E_{8,a}=-a+126a^2-1260a^3+2520a^4,
\\&E_{10,a}=-a+510a^2-13230a^3+75600a^4-113400a^5\endaligned\tag 1.2$$
\par In Section 2 we deduce some identities and an inversion formula
involving $\{E_{n,a}\}$. In Sections 3-5 we establish congruences
for $E_{2n,a}\mod{2^{{\rm ord}_2n+8}}$, $E_{2n,a}\mod{3^{{\rm
ord}_3n+5}}$ and $E_{2n,a}\mod{5^{{\rm ord}_5n+4}}$ provided that
$a$ is a nonzero integer, where ${\rm ord}_pn$ is the least
nonnegative integer $\alpha$ such that $p^{\a}\mid n$ but
$p^{\a+1}\nmid n$. See Theorems 3.1, 4.1, 4.2 and 5.1.
\par In addition to the above notation, throughout this paper we use
$\Bbb Z$ and $\Bbb N$ to denote the set of integers and the set of
positive integers, respectively.
\section*{2. Identities involving $\{E_{n,a}(x)\}$}
\pro{Definition 2.1} For any given real number $a$ we define
$\{E_{n,a}\}$ and $\{E_{n,a}(x)\}$ by
$$E_{0,a}=1, \q
E_{n,a}=-a\sum_{k=1}^{[n/2]}\binom n{2k}E_{n-2k,a}\ (n\ge 1)$$
 and $$E_{n,a}(x)=\s_{k=0}^n\b
nkE_{k,a}x^{n-k}\ (n\ge 0).$$ \endpro \par By the definition, we
have
$$E_{2n-1,a}=0\qtq{and}\sum_{k=1}^n\b{2n}{2k}E_{2n-2k,a}=-\f
1aE_{2n,a}\qtq{for}a\not=0.\tag 2.1$$

\pro{Theorem 2.1} For any real number $a$ we have
$$\sum_{n=0}^{\infty}E_{n,a}\frac{t^n}{n!}
=\frac1{\frac a2(e^t+e^{-t})+1-a}\q \t{and}\q
\s_{n=0}^{\i}E_{n,a}(x)\frac {t^n}{n!}=\frac {e^{xt}}{\frac
a2(e^t+e^{-t})+1-a}.$$\endpro Proof.
 Since
$$\aligned&\Big(\frac a2(e^t+e^{-t})+1-a\Big)\Big(\sum_{m=0}^{\infty}E_{m,a}\frac{t^m}{m!}\Big)
\\&=\Big(1+a\sum_{k=1}^{\infty}\frac {t^{2k}}{(2k)!}\Big)\Big(\sum_{m=0}^{\infty}E_{m,a}\frac{t^m}{m!}\Big)
\\&=1+\sum_{n=1}^{\infty}\Big(a\sum_{k=1}^{[n/2]}
\frac1{(2k)!}\cdot \f{E_{n-2k,a}}{(n-2k)!}+\f{E_{n,a}}{n!}\Big)t^n
\\&=1+\sum_{n=1}^{\infty}\Big(a\sum_{k=1}^{[n/2]}\binom n{2k}E_{n-2k,a}+E_{n,a}\Big)\frac{t^n}{n!}
=1,\endaligned$$ we deduce the first result.  To complete the proof,
we note that
 $$\s_{n=0}^{\i}E_{n,a}(x)\frac {t^n}{n!}
=\Big(\s_{k=0}^{\i}E_{k,a}\frac
{t^k}{k!}\Big)\Big(\s_{m=0}^{\i}x^m\frac {t^m}{m!}\Big)=\frac1{\frac
a2(e^t+e^{-t})+1-a}\cdot \t{e}^{tx}.$$ \pro{Corollary 2.1} Let $a$
be a real number. Then
$$E_{n,a}(1-x)=\sum_{k=0}^n\b nk(-1)^{k}E_{k,a}(x).$$
\endpro
Proof. By Theorem 2.1,
$$\aligned\sum_{n=0}^{\infty}E_{n,a}(1-x)\f{(-t)^n}{n!}
&=\f{e^{(1-x)(-t)}}{\f a2(e^t+e^{-t})+1-a}=e^{-t}\f{e^{xt}}{\f
a2(e^t+e^{-t})+1-a}
\\&=\Big(\sum_{k=0}^{\infty}(-1)^k\f{t^k}{k!}\Big)\Big(\sum_{k=0}^{\infty}E_{k,a}(x)\f{t^k}{k!}\Big)
\\&=\sum_{n=0}^{\infty}\Big(\sum_{k=0}^n\b nk(-1)^{n-k}E_{k,a}(x)\Big)\f{t^n}{n!}
\endaligned$$
Comparing the coefficients of $t^n$ on both sides we deduce the
result. \pro{Theorem 2.2} For any real number $a$ and positive
integer $n$ we have
$$\align&({\rm i})\q\f a2\big(E_{n,a}(x+1)+E_{n,a}(x-1)\big)+(1-a)E_{n,a}(x)=x^n,
\\&({\rm ii})\q\sum_{k=0}^n\binom nkE_{n-k,a}(x)\left\{\frac a
2((y+1)^k+(y-1)^k)+(1-a)y^k\right\}=(x+y)^n,
\\&({\rm iii})\q E_{n,a}(x)=x^n-a\s_{k=1}^{[n/2]}\b n{2k}E_{n-2k,a}(x).
 \endalign$$
\endpro
Proof. By Theorem 2.1,
 $$\aligned&\sum_{n=0}^{\infty}\Big(\f
 a2\big(E_{n,a}(x+1)+E_{n,a}(x-1)\big)+(1-a)E_{n,a}(x)\Big)\f{t^n}{n!}
 \\&\q=\f{\f a2\big(\te^{(x+1)t}+\te^{(x-1)t}\big)+(1-a)\te^{tx}}
 {\f a2(\te^t+\te^{-t})+1-a}
 =\te^{xt}=\sum_{n=0}^{\infty}x^n\f{t^n}{n!}.\endaligned$$
 Comparing the coefficients of $x^n$ on both sides yields (i).
 \par Using Definition 2.1 and Theorem 2.1 we see that
$$\align \sum_{n=0}^{\infty}E_n(x+y)\f{t^n}{n!}
&=\f{\te^{(x+y)t}}{\f a2(\te^t+\te^{-t})+1-a}
=\Big(\sum_{k=0}^{\infty}y^k\f{t^k}{k!}\Big)\Big(\sum_{m=0}^{\infty}
E_m(x)\f{t^m}{m!}\Big) \\&=\sum_{n=0}^{\infty}\Big(\sum_{k=0}^n\b
nkE_{n-k}(x)y^k\Big)\f{t^n}{n!}.\endalign$$ Thus,
$$E_n(x+y)=\sum_{k=0}^n\b
nkE_{n-k}(x)y^k=\sum_{k=0}^n\b nkE_k(x)y^{n-k}.\tag 2.2$$ This
together with (i) yields (ii). Taking $y=0$ in (ii) and applying the
fact $E_{2k-1,a}=0$ we deduce (iii). The proof is now complete.

\par{\bf Remark 2.1.} In the case $a=2$, Theorem 2.2(i) is known.
See [8, p.427]. \pro{Corollary 2.2} For any real number $a$ and
positive integer $n$ we have
$$\align &\f a2\sum_{k=1}^n\b{2n}{2k}E_{2n-2k,a}
((x+1)^{2k}+(x-1)^{2k})+(1-a)\sum_{k=1}^n\b{2n}{2k}E_{2n-2k,a}
x^{2k}\\&=x^{2n}-E_{2n,a}.\endalign$$
\endpro
Proof. Taking $x=0$ and then substituting $n,y$ with $2n,x$ in
Theorem 2.2(ii) we obtain
$$\f a2\sum_{k=0}^n\b{2n}{2k}E_{2n-2k,a}
((x+1)^{2k}+(x-1)^{2k})
+(1-a)\sum_{k=0}^n\b{2n}{2k}E_{2n-2k,a}x^{2k}=x^{2n}.$$ This yields
the result.

 \pro{Theorem 2.3} Let $a$ be a real number with $a\not=0$
  and $n\in\Bbb N$.  Then
  $$\sum_{k=1}^n\b {2n}{2k} 2^{2k}E_{2n-2k,a}=\f 2a+\f {2-4a}{a^2}
  E_{2n,a}$$
  and so
  $$E_{2n,a}=\frac a {2a-1}-\frac {a^2n}{2(2a-1)}\sum_{k=1}^n\binom{2n-1}{2k-1}
\frac {4^k}kE_{2n-2k,a}\qtq{for}a\not=\f 12.$$
  \endpro
  Proof. Taking $x=1$ in Corollary 2.2 and then
  applying (2.1) and the fact that $\b{2n}{2k}=\f
nk\b{2n-1}{2k-1}$   we  deduce the result.

  \pro{Theorem 2.4} Let $a$ be a real number with $a\not=0$
   and  $n\in\Bbb N$.  Then
  $$\sum_{k=1}^n\b {2n}{2k} 3^{2k}E_{2n-2k,a}
  =\f {2^{2n+1}}a+\f{4(a-1)}{a^2}-\f{(3a-2)^2}{a^3}E_{2n,a}.$$
  \endpro
  Proof. Taking $x=2$ in Corollary 2.2 and then
   applying (2.1) and Theorem 2.3 we deduce the result.

\pro{Theorem 2.5} Let $a$ be a real number with $a\not=0$
 and $n\in\Bbb N$.  Then
  $$\align &\sum_{k=1}^n\b {2n}{2k} 4^{2k}E_{2n-2k,a}
  \\&=\f {2(3^{2n}-1)}a+\f{8(a-1)^2}{a^3}+\f{2^{2n+2}(a-1)}{a^2}
  -8\f{(a-1)^2(2a-1)}{a^4}E_{2n,a}.\endalign$$
  \endpro
  Proof. Taking $x=3$ in Corollary 2.2  and then applying
    Theorems 2.3 and 2.4 we deduce the result.

\pro{Theorem 2.6} Let $a$ be a real number with $a\not=0$
 and $n\in\Bbb N$.  Then
  $$\align \sum_{k=1}^n\b {2n}{2k}5^{2k} E_{2n-2k,a}
  &=\f {2\cdot
  4^{2n}}a+\f{(3a^2-8a+4)2^{2n+1}}{a^3}+\f{4(a-1)3^{2n}}{a^2}
  \\&\q+\f{8(a-1)(a^2-4a+2)}{a^4}-\f{(5a^2-10a+4)^2}{a^5}E_{2n,a}.\endalign$$
  \endpro
  Proof. Taking $x=4$ in Corollary 2.2
   and then applying  Theorems 2.4 and 2.5 we deduce the result.

 \pro{Theorem 2.7} Let $a$ be a real number. For two sequences $\{a_n\}$ and $\{b_n\}$ we have the following
 inversion formula:
 $$\align &b_n=\s_{k=0}^n\b nk\Big((1-a)(-x)^k+\frac a2\big((1-x)^k+(-1-x)^k\big)\Big)a_{n-k}
\\& \iff a_n=\s_{k=0}^n \b nkE_{k,a}(x)b_{n-k}.\endalign $$
\endpro
 Proof. Clearly
 $$\align &\frac {1-a+\frac a2(e^t+e^{-t})}{e^{xt}}\Big(\s_{n=0}^{\i}a_n\frac {t^n}{n!}\Big)
 \\&\q=\Big((1-a)e^{-xt}+\frac a2 (e^{(1-x)t}+e^{(-1-x)t})\Big)
 \Big(\s_{m=0}^{\i}a_m\frac
 {t^m}{m!}\Big)
 \\&\q=\s_{k=0}^{\i}\frac {t^k}{k!}\Big((1-a)(-x)^k)
 +\frac a2\big((1-x)^k+(-1-x)^k\big)\Big)\Big(\s_{m=0}^{\i}a_m
 \frac{t^m}{m!} \Big )
 \\&\q=\s_{n=0}^{\i}\Big(\s_{k=0}^n\b nk\big((1-a)(-x)^k+\frac
 a2\big((1-x)^k+(-1-x)^k\big)a_{n-k}\Big)\frac {t^n}{n!}.
 \endalign $$
Thus, using Theorem 2.1 we see that
$$\align &b_n=\s_{k=0}^n\b nk\Big((1-a)(-x)^k+\frac a2\big((1-x)^k+(-1-x)^k\big)
\Big)a_{n-k}
\\&\Leftrightarrow
\frac {1-a+\frac a2(e^t+e^{-t})}{e^{xt}}\Big(\s_{n=0}^{\i}a_n\frac
{t^n}{n!}\Big)=\s_{n=0}^{\i}b_n\frac {t^n}{n!}
\\&\Leftrightarrow \s_{n=0}^{\i}a_n\frac
{t^n}{n!}=\Big(\s_{n=0}^{\i}b_n\frac
{t^n}{n!}\Big)\Big(\s_{n=0}^{\i}E_{n,a}(x)\frac {t^n}{n!}\Big)
\\&\Leftrightarrow a_n=\s_{k=0}^n \b nkE_{k,a}(x)b_{n-k}.\endalign $$
This proves the theorem.

\pro{Corollary 2.3} Let $a$ be a real number. For two sequences
$\{a_n\}$ and $\{b_n\}$ we have the following
 inversion formula:
 $$\align &b_n=a\sum_{k=0}^{[n/2]}\b n{2k}a_{n-2k}+(1-a)a_n\
 (n=0,1,2,\ldots)
\\& \iff a_n=\sum_{k=0}^{[n/2]}\b n{2k}E_{2k,a}b_{n-2k}\
 (n=0,1,2,\ldots).\endalign$$
 \endpro
 Proof. Taking $x=0$ in Theorem 2.7 we derive the result.
 \par\q
 \par{\bf Remark 2.2.} In the case $a=2$, Corollary 2.3 was given
 by the first author in [7].
\pro{Theorem 2.8} Let $a$ be a nonzero real number and $n\in\Bbb N$.
Then
$$\align &\sum_{k=1}^{[n/4]}\b n{4k}E_{n-4k,a}(x)(2(1-a)+a(-4)^k)
\\&=\f 12((x+i)^n+(x-i)^n)+\f{1-a}ax^n-\f 1aE_{n,a}(x).\endalign$$
\endpro
Proof. Taking $y=\pm i$ in Theorem 2.2(ii) we find
$$\sum_{k=0}^n\b nkE_{n-k,a}(x)\Big\{\f a2((\pm i+1)^k+(\pm
i-1)^k)+(1-a)(\pm i)^k\Big\}=(x\pm i)^n.$$ Thus,
$$\align &\sum_{k=0}^n\b nkE_{n-k,a}(x)\Big\{\f
a2((i+1)^k+(i-1)^k+(-i+1)^k+(-i-1)^k)\\&\q+(1-a)(i^k+(-i)^k)\Big\}
=(x+i)^n+(x-i)^n.\endalign$$ Observe that
$$(i+1)^k+(i-1)^k+(-i+1)^k+(-i-1)^k=\cases (-1)^{\f k4}2^{\f k2+2}
&\t{if $4\mid k$,}\\0&\t{if $4\nmid k$}
\endcases$$ and
$$i^k+(-i)^k=\cases 2(-1)^{\f k2}&\t{if $2\mid k$,}
\\ 0&\t{if $2\nmid k$.}\endcases$$
We then have
$$\align &2a\sum_{k=0}^{[n/4]}\b n{4k}E_{n-4k,a}(x)(-1)^k2^{2k}
+2(1-a)\sum_{k=0}^{[n/2]}\b n{2k}E_{n-2k,a}(x)(-1)^k
\\&=(x+i)^n+(x-i)^n.\endalign$$
By Theorem 2.2(iii),
$$\align &a\sum_{k=0}^{[n/2]}\b n{2k}E_{n-2k,a}(x)(-1)^k
+x^n-(1-a)E_{n,a}(x)
\\&=a\sum_{k=0}^{[n/2]}\b n{2k}E_{n-2k,a}(x)((-1)^k+1)
=2a\sum_{k=0}^{[n/4]}\b n{4k}E_{n-4k,a}(x).\endalign$$ Thus,
$$\align &2a\sum_{k=0}^{[n/4]}\b n{4k}E_{n-4k,a}(x)(-4)^k
+2(1-a)\Big(2\sum_{k=0}^{[n/4]}\b
n{4k}E_{n-4k,a}(x)-\f{x^n}a+\f{1-a}aE_{n,a}(x)\Big)\\&=(x+i)^n+(x-i)^n.
\endalign$$
This yields the result.
 \pro{Corollary 2.4} Let $a$ be a real number with $a\not=0$ and $n\in\Bbb
N$. Then
$$\sum_{k=1}^{[n/2]}\binom {2n}{4k}E_{2n-4k,
a}\big(a(-4)^k+2(1-a)\big)=(-1)^n-\f 1aE_{2n,a}.$$
\endpro
Proof. Taking $x=0$ and replacing $n$ with $2n$ in Theorem 2.8 we
deduce the result.
\par\q
\par{\bf Remark 2.3.} In the case $a=1$, Corollary 2.4 can be
found in [2, p.643]. In the case $a=2$, Corollary 2.4 was proved by
the first author in [7, Theorem 2.4(iv)].

\section*{3. A congruence for $E_{2n,a}\mod{2^{{\rm
ord}_2n+8}}$}
\par\q Let $a$ be a nonzero integer and $n\in\Bbb N$. By Theorem 2.3,
$$E_{2n,a}=\frac a {2a-1}-\frac {a^2n}{2(2a-1)}\sum_{k=1}^n\binom{2n-1}{2k-1}
\frac {4^k}kE_{2n-2k,a}.\tag 3.1$$ Since $E_{2n,a}\in\Bbb Z$, we get
$$E_{2n,a}\e a-\f{a^2n}2\cdot (2n-1)4E_{2n-2,a}\e a\mod 2.\tag 3.2$$
From (3.1) and (3.2) we see that
$$E_{2n,a}\e \f a{2a-1}-\f{a^2n}{2(2a-1)}(2n-1)\cdot 4a
\e a(1-2n)\mod 4.\tag 3.3$$ As $E_{2m,a}\in\Bbb Z$ and
$\f{4^k}{2k}\e 0\mod{2^5}$ for $k\ge 3$, from (3.1) we see that for
$n\ge 2$,
$$\aligned &E_{2n,a}-\frac a {2a-1}\\&\e -\frac
{a^2n}{2(2a-1)}\Big(\b{2n-1}14E_{2n-2,a}+\b{2n-1}38E_{2n-4,a}\Big)
\\&= -\f{a^2}{2a-1}\Big(2n(2n-1)E_{2n-2,a}
+\f{n(n-1)(2n-1)(2n-3)}68E_{2n-4,a}\Big)\mod {32}.\endaligned\tag
3.4$$
 Thus,
$$\align E_{2n,a}&
\e \f a{2a-1}-\f{a^2}{2a-1}n(2n-1)2a(1-2(n-1))
\\&=2a^2-a+2a^3(2a-1)n(2n-1)(2n-3)
\e 2a^2-a-2a^3(2a-1)n\\&\e 2a^2-a-(4a^3-2a^3)n \mod 8.\endalign$$
That is,
$$E_{2n,a}\e 2a^2-a-2a^3n \mod 8\qtq{for} n\ge 2.\tag 3.5$$
From (3.2), (3.4) and (3.5) we see that for $n\ge 3$,
$$\align &E_{2n,a}-\f a{2a-1}\\&\e
-\f{a^2}{2a-1}\Big(2n(2n-1)(2a^2-a-2a^3(n-1))
+\f{n(n-1)(2n-1)(2n-3)}68a\Big)
\\&\e -2n(2n-1)a^3+\f{2a^5}{2a-1}\cdot
2n(n-1)(2n-1)-\f{a^2}{2a-1}\cdot\f{n(n-1)}2\cdot 8a
\\&\e -2n(2n-1)a^3\mod{16}.\endalign$$
This is also true for $n=2$ by (1.2). Thus,
$$E_{2n,a}\e  \f a{2a-1}-2n(2n-1)a^3\mod{16}\qtq{for}n\ge 2.\tag 3.6$$
From (3.2), (3.4) and (3.6) we see that for $n\ge 3$,
$$ \align E_{2n,a}-\f
a{2a-1}&\e -\f{a^2}{2a-1}\Big(2n(2n-1)\Big(\f
a{2a-1}-(2n-2)(2n-3)a^3\Big)\\&\q+\f{n(n-1)(2n-1)(2n-3)}68a(1-2(n-2))\Big)
\\&\e-\f{a^3}{(2a-1)^2}2n(2n-1)-8a^3\cdot
\f{n(n-1)}2-\f{8a^3}{2a-1}(1-2n)\cdot\f{n(n-1)}2
\\&\e -\f{a^3}{(2a-1)^2}2n(2n-1)-8a^3\cdot
\f{n(n-1)}2-8a^3(1-2n)\cdot\f{n(n-1)}2
\\&\e -a^3(1-4a(a-1))2n(2n-1)-8a^3\cdot \f{n(n-1)}2(1+1-2n)
\\&\e -2a^3n(2n-1)+8a^2n-8a^3n+8a^3n(n-1)^2\\&= 2a^3n(4n^2-10n+1)+8a^2n
\\&\e 8a^3n^3+12a^3n^2+(2a^3+8a^2)n \mod{32}.\endalign$$  As
$$E_{4,a}=6a^2-a\e \f a{2a-1}+64a^3+48a^3+2(2a^3+8a^2)\mod{32},$$ we
obtain
$$E_{2n,a}\e \f
a{2a-1}+8a^3n^3+12a^3n^2+(2a^3+8a^2)n\mod{32}\qtq{for}n\ge 2.\tag
3.7$$ From (3.1) we see that for $n\ge 4$,
$$ E_{2n,a}\e \f
a{2a-1}+S_1+S_2+S_3+S_4\mod{2^{{\rm ord}_2n+8}},\tag 3.8$$ where
 $$\aligned &S_1=-\f{a^2n}{2a-1}2(2n-1)E_{2n-2,a},
 \\&S_2=-\f{a^2}{2a-1}\cdot\f{n
(n-1)(2n-1)(2n-3)}34E_{2n-4,a},
\\&S_3=-\f{a^2}{2a-1}\cdot\f{n
(n-1)(n-2)(2n-1)(2n-3)(2n-5)}{45}16E_{2n-6,a},
\\&\ S_4=-\f{a^2}{2a-1}\cdot\f{n
(n-1)(n-2)(n-3)(2n-1)(2n-3)(2n-5)(2n-7)}{7\cdot5\cdot9}16E_{2n-8,a}.
\endaligned\tag 3.9$$
Thus,  using (3.2)-(3.8) we see that for $n\ge 4$,
$$ \align &E_{2n,a}-\f
a{2a-1}\\&\e-\f{a^2}{2a-1}\Big(2n(2n-1)(\f
a{2a-1}+2a^3(n-1)(4(n-1)^2-10(n-1)+1)+8a^2(n-1))
\\&\q+\f{n(n-1)(2n-1)(2n-3)}68(2a^2-a-2a^3(n-2))
\\&\q+\f{n(n-1)(n-2)(2n-1)(2n-3)(2n-5)}{5\cdot9\cdot2}32a\Big)
\\&\e-\f{a^3}{(2a-1)^2}2n(2n-1)-\f{n(n-1)}2\Big\{\f{8a^5}{2a-1}(2n-1)(4(n-1)^2-2(n-1)+1+32a)
\\&\q-\f{a^2}{2a-1}(-32an+16a^2-8a-16a^3n+32a)-32a(n-2)\Big\}
\\&\e-\f{a^3}{(2a-1)^2}2n(2n-1)+\f{n(n-1)}2\Big\{\f1{2a-1}(-16a^5n+32a(n-1)
-16a^5(n-1)\\&\q+8a^5+32a) +(32a-8a^3-16a^5n)-32a(n-2)\Big\}
\\&\e-\f{a^3}{(2a-1)^2}2n(2n-1)+(48a^6-24a^5+32a-8a^3-16a^5n-32an)\f{n(n-1)}2
\\&\e-\f{a^3}{(2a-1)^2}2n(2n-1)+(48a^2-24a^3+32a-8a^3-16a^3n-32an)\f{n(n-1)}2
\\&\e-a^3(1-4a(a-1))2n(2n-1)+(-16a^2-16a^3n-32an)\f{n(n-1)}2
\\&\e 2a^3(4a^2-4a-1)n(2n-1)+16a^2(an-1)\f{n(n-1)}2
\\&=8a^3n^3+(16a^5-16a^4-12a^3-8a^2)n^2+(-8a^5+8a^4+2a^3+8a^2)n\mod{64}.\endalign$$
and so
$$E_{2n,a}\e \f
a{2a-1}+8a^3n^3+(4a^3-24a^2)n^2+(8a^4-6a^3+8a^2)n\mod{64}\qtq{for}n\ge
4.\tag 3.10$$  As
$$\align S_1&\e-\f{a^2}{2a-1}2n(2n-1)\Big(\f a{2a-1}+8a^3(n-1)^3+(4a^3-24a^2)(n-1)^2
\\&\q+(8a^4-6a^3+8a^2)(n-1)\Big)
\\&\e-\f{a^3}{(2a-1)^2}2n(2n-1)-\f a{2a-1}\f{n(n-1)}2(32a^3(2n-1)(n-1)^2
\\&\q+16a^5(2n-1)(n-1)
 +32a^2(n-1)(2n-1)+64a(2n-1)-24a^5n(2n-1))
\\&\e-a^3(1-4a(a-1))2n(2n-1)\\&\q-(-32a^3-32a^3n+40a^5-32a^2n+32a^2-64a)
\f1{2a-1}\cdot\f{n(n-1)}2
\\&\e-a^3(1-4a(a-1))2n(2n-1)-(-32a^3(2a-1)-32a^3n(2a-1)\\&\q-40a^5(2a+1)(4a^2+1)
-32a^2n(2a-1)+32a^2(2a-1)-64a)\f{n(n-1)}2
\\&\e-4a^3n^2+2a^3n+16a^5n^2-8a^5n-16a^4n^2+8a^4n+(-64a+32a^3-64an+32a^3n\\&\q-64a
-80a^6-32a^3-40a^5-64an+32a^2n+64a-32a^2-64a)\f{n(n-1)}2
\\&\e-4a^3n^2+2a^3n+16a^5n^2-8a^5n-16a^4n^2+8a^4n+(-32a^3n-48a^6+40a^5
\\&\q-32a^2n+32a^2)\f{n(n-1)}2\mod{128},\endalign$$
$$\align S_2&\e-\f{a^2}{2a-1}\f{n (n-1)(2n-1)(2n-3)}68\Big(\f a{2a-1}-2(n-2)(2(n-2)-1)a^3\Big)
\\&\e-\f{a^3}{(2a-1)^2}\f{8(2n-1)(2n-3)}3\f{n(n-1)}2-\f{a^5}{3(2a-1)}(-32n^2
+16n-32)\f{n(n-1)}2
\\&\e(32a^3n^2+64an+24a^3-32a^2)\f{n(n-1)}2+(32a^3n^2+32a^3+48a^5n)\f1{2a-1}
\\&\e(32a^3n^2+64an+24a^3-32a^2)\f{n(n-1)}2+(32a^3n^2(2a-1)+32a^3(2a-1)
\\&\q-48a^5n(2a+1)(4a^2+1))\f{n(n-1)}2
\\&\e(-8a^3-32a^2+64a+32a^2n-64an-48a^5n)\f{n(n-1)}2\mod{128},\endalign$$
$$\align S_3&\e-\f{a^2}{2a-1}\f{n
(n-1)(n-2)(2n-1)(2n-3)(2n-5)}{5\cdot9\cdot2}32a(1-2(n-3))
\\&\e-\f{a^3}{2a-1}32(n-2)(2n-3)(1-2(n-3))\f{n(n-1)}2
\\&\e-32a^3(2a-1)(n-2)(2n-3)(1-2(n-3))\f{n(n-1)}2
\\&\e(-64an+32a^3n-64a)(2n-3)(1-2(n-3))\f{n(n-1)}2
\\&\e(32a^3n+64a)(1-2(n-3))\f{n(n-1)}2
\e(64a+32a^3n)\f{n(n-1)}2\mod{128}\endalign$$ and
$$\align S_4&\e-\f{a^2}{2a-1}\f{n
(n-1)(n-2)(n-3)(2n-1)(2n-3)(2n-5)(2n-7)}{7\cdot5\cdot9\cdot2}32a
\\&\e-32a(n-2)(n-3)\f{n(n-1)}2
\e(-32an^2+32an-64a)\f{n(n-1)}2\mod{128},\endalign$$ using (3.8) we
deduce that
$$\aligned E_{2n,a}&\e\f
a{2a-1}-16an^4-24a^5n^3-(24a^6-60a^5+16a^4+8a^3-48a)n^2
\\&\q+(24a^6-28a^5+8a^4+6a^3-32a)n\mod{128}\qtq{for}n\ge 4.
\endaligned\tag 3.11$$

\pro{Theorem 3.1} Let $a$ be a nonzero integer, $n\in\Bbb N$ and
$n\ge 5$. Then
$$\align E_{2n,a}&\e \f a{2a-1}-96a^3n^5+(16a^5-32a^4-64a^2)n^4+(72a^7-64a^3)n^3
\\&\q-(24a^7-120a^6+92a^5-56a^3+128a)n^2
\\&\q-(80a^7-72a^6-20a^5+104a^4-6a^3+64a^2-128a)n\mod{2^{{\rm
ord}_2n+8}}.\endalign$$
\endpro
Proof. Set $\alpha={\rm ord}_2n.$ Let $S_1,S_2,S_3$ and $S_4$ be
given by (3.9).  Since
$$\align \f 1{(2a-1)^2}&=\f{1-4a(a-1)}{1-16a^2(a-1)^2}
\e (1-4a(a-1))(1+16a^2(a-1)^2)
\\&\e (1-4a(a-1))(1+32a(a-1))\e 1+28a(a-1)
\mod{128}\endalign$$ we see that

 $$\align S_1&\e-\f{a^2n}{2a-1}2(2n-1)(\f a{2a-1}-16a(n-1)^4-24a^5(n-1)^3
 \\&\q-(24a^6-60a^5+16a^4+8a^3-48a)(n-1)^2
 \\&\q+(24a^6-28a^5+8a^4+6a^3-32a)(n-1))
  \\&\e-a^32n(2n-1)(1+28a(a-1))+a^2n(2a+1)
  (4a^2+1)(16a^4+1)(32a^5+128a\\&\q+96a^6+48a^4
 +28a^3+32an^4+192an^3+96an^2+128an+64a^5n^3+136a^5n^2
 \\&\q+120a^5n+80a^6n^2+176a^6n
 +192a^4n^2+80a^4n+104a^3n^2+156a^3n+192an^5\\&\q+160a^5n^4+160a^6n^3
 +192a^4n^3+224a^3n^3)
 \\&\e130a^3n-56a^4n-16a^7n+104a^6n+128an-96n^3a^3+64a^2n+44a^5n^2+112a^4n^2
 \\&\q-120a^7n^2-120a^6n^2-32n^4a^3-112a^6n^3-88a^7n^3+32a^4n^4-64n^5a^3
 \\&\q+124a^3n^2-64a^2n^2+64a^2n^3-24a^5n^3-96a^4n^3+84a^5n\mod{2^{\alpha+8}},\endalign$$
$$\align S_2&\e-\f{a^2n}{2a-1}\cdot\f{4(n-1)(2n-1)(2n-3)}3\Big(\f
a{2a-1}+8a^3(n-2)^3
\\&\q+(4a^3-24a^2)(n-2)^2+(8a^4-6a^3+8a^2)(n-2)\Big)
\\&\e-a^3n(1+4a(1-a))\f{4(n-1)(2n-1)(2n-3)}3\\&\q+4a^2n(2a+1)(4a^2+1)(16a^4+1)
\f{ (n-1)(2n-1)(2n-3)}3\\&\q\times(8a^3(n-2)^3+(4a^3-24a^2)(n-2)^2
+(8a^4-6a^3+8a^2)(n-2))
\\&\e -60a^3n+48a^4n+128an+112n^3a^3-104a^5n^2+16a^4n^2-80a^6n^2+80a^5n^4
\\&\q-48a^6n^3-48n^4a^3-4a^3n^2-64a^2n^2+128an^2-104a^5n^3\mod{2^{\alpha+8}},\endalign$$
$$\align S_3&\e-\f{ a^2(2a-1)(n-1)(n-2)(2n-1)(2n-3)(2n-5)}516(2a^2-a-2a^3(n-3))
\\&\e16a^3n^2-96a^3n-112n^3a^3-64a^2n^3+64a^2n^4-64n^4a^3\mod{2^{\alpha+8}}\endalign$$
and
$$\align S_4&\e\f {16}5a^2n(2a-1)(n-1)(n-2)(n-3)(2n-1)(2n-3)(2n-5)(2n-7)
\\&\q\times(2a^2-a-2a^3(n-4))
\\&\e112a^3n^2+64a^2n^2-64a^2n^4-48n^4a^3-96n^5a^3+32a^3n\mod{2^{\alpha+8}}\endalign$$ Now combining
the above with (3.8) yields the result. \pro{Corollary 3.1} Let $a$
be a nonzero integer and $n\in\Bbb N$ with $n\ge 5$.
\par $(\t{\rm i})$ If $2\mid a$, then
$$E_{2n,a}\e\f
a{2a-1}+(4a^5-8a^3)n^2+(4a^5-8a^4+6a^3)n\mod{2^{{\rm ord}_2n+8}}.$$
 \par $(\t{\rm ii})$
If $2\nmid a$, then
$$\align E_{2n,a}&\e\f
a{2a-1}-96(a+a^2-1)n^5+(16a-96a^2)n^4\\&\q+(-104a^3+112a)n^3
+(56a^3-116a+104a^2-112)n^2\\&\q+(62a^3-116a-56a^2+88)n\mod{2^{{\rm
ord}_2n+8}}.\endalign$$
\endpro
\section*{4. A congruence for $E_{2n,a}\mod{3^{{\rm ord}_3n+5}}$}
 \par\q Let $a$ be a nonzero integer and $n\in\Bbb N$. By Theorem
2.4, we have
$$ E_{2n,a}=\frac 1{(3a-2)^2}\Big(2^{2n+1}a^2+4a(a-1)
-a^3n\sum_{k=1}^n\b{2n-1}{2k-1}\f{3^{2k}}kE_{2n-2k,a}\Big).\tag
4.1$$ Since $E_{2n,a}\in\Bbb Z$ and $2^{2n}=(1+3)^n\e 1+3n\mod 9$,
using (4.1)we see that for $n\ge 2$,
$$\align E_{2n,a}&\e\f 1{2(1-3a)}((1+3n)a^2+2a^2-2a)
\\&\e-4(1+3a)(3a^2-2a+3na^2)\e3a^2-a-3a^2n\mod 9.\tag 4.2
\endalign$$
 \pro{Theorem 4.1} Let $a$ be a nonzero integer with $3\mid a$, $n\in\Bbb N$
and $n\ge 2$. Then
$$\align E_{2n,a}&\e \f{2a}{3a-2}+ 9a^2n^3+9a^2n^2-3a^2n
\mod{3^{{\rm ord}_3n+5}}.\endalign$$
\endpro
Proof. It is clear that
  $$\frac 1{(3a-2)^2}=\frac
{(3a+2)^2}{(9a^2-4)^2}\e -5(9a^2+12a+4)\e 21a-20 \mod{81}$$ and
$$2^{2n}=(1+3)^n=1+n\sum_{k=1}^n\b {n-1}{k-1}\f{3^k}k.$$
 Thus,
$$\align E_{2n,a}&\e\frac
1{(3a-2)^2}\Big(2a^2(1+n\sum_{k=1}^n\b{n-1}{k-1}\f{3^k}k)+4a(a-1)\Big)
\\&\e\frac 1{(3a-2)^2}\Big(6a^2-4a+2a^2n\big(3+\f 92(n-1)+\f 92(n-1)(n-2)\big)\Big)
\\&\e\frac {6a^2-4a}{(3a-2)^2}+n\big(6a^2\cdot\f 1{4-3a}+9a^2(n-1)+9a^2(n-1)(n-2)\big)
\\&\e \f{2a}{3a-2}+n(6a^2(7+3a)+9a^2+9a^2n^2-18a^2n)
\\&\e \f{2a}{3a-2}+9a^2n^3+9a^2n^2-3a^2n
\mod{3^{{\rm ord}_3n+5}}.
\endalign$$
This proves the theorem.

 \pro{Theorem 4.2} Let $a$ be a nonzero integer with $3\nmid a$, $n\in\Bbb N$ and
$n\ge 3$. Then
$$\align E_{2n,a}&\e\f{2a}{3a-2}+(54a^3-99a^2+27a-81)n^3-(9a^4-27a^3-63a^2
\\&\q-81a-108)n^2-(117a^4+117a^3+111a^2-54a-108)n\mod{3^{{\rm
ord}_3n+5}}.\endalign$$
\endpro
Proof. Clearly $(3a-2)^2(7-6a)=135a^2-108a+28-54a^3\e 1\mod{27}$ and
$$\align &(3a-2)^2(-27a^3+108a^2-60a-20)\\&=-243a^5+1296a^4-1944a^3+972a^2-80\e
81a^4+163\e 1\mod {3^5}.\endalign$$
 Now, using (4.1), (4.2) and the above we deduce that
$$\align E_{2n,a}&\e(7-6a)\Big(2\big(1+3n+\f{9n(n-1)}2\big)a^2+4a(a-1)-9an(2n-1)E_{2n-2,a}\Big)
\\&\e(7-6a)\Big((2+6n+9n(n-1))a^2+4a^2-4a+9n(2n-1)\Big)
\\&\e(7-6a)(6a^2+6a^2n+9n-4a)
\\&\e12a^2+15a^2n+9n+26a+18a^3+18a^3n
\\&\e(-12a^2-9a+9)n+12a^2-10a
\mod{27}\endalign$$  and so
$$\align E_{2n,a}&\e\frac
1{(3a-2)^2}\Big(2a^2\Big(1+n\sum_{k=1}^n\b{n-1}{k-1}\f{3^k}k\Big)
+4a(a-1)\\&\q-a^3n\sum_{k=1}^n\b{2n-1}{2k-1}\f{3^{2k}}kE_{2n-2k,a}\Big)
\\&\e\frac 1{(3a-2)^2}\Big(6a^2-4a+2a^2n\big(3+\f 92(n-1)+\f 92(n-1)(n-2)
\\&\q+\f {27}8(n-1)(n-2)(n-3)\big)
\\&\q-a^3n\Big(9(2n-1)E_{2n-2,a}+\f{27}2(n-1)(2n-1)(2n-3)E_{2n-4,a}\Big)\Big)
\\&\e\frac 1{(3a-2)^2}\Big(6a^2-4a+2a^2n\big(3+\f 92(n-1)+\f 92(n-1)(n-2)
\\&\q+\f {27}8(n-1)(n-2)(n-3)\big)
\\&\q-a^3n\Big(9(2n-1)\big(-12a^2-9a+9)(n-1)+12a^2-10a\big)
\\&\q+\f{27}2(n-1)(2n-1)(2n-3)(3a^2-a-3a^2(n-2))\Big)
\\&\e \f{2a}{3a-2}+(-27a^3+108a^2-60a-20)(96na^2+117a^2n^2+90a^2n^3
\\&\q+189n^4a^2+54n^3a^5
+189n^2a^5+207a^4n^2+81a^3n^3
\\&\q+216na^5+72a^4n+162na^3+54n^4a^4+162n^5a^5)
\\&\e\f{2a}{3a-2}+24a^2n+90a^2n^2+81an^2-9a^4n^2
-108n^3a^5\\&\q-54na^5-63a^4n -99a^2n^3+81n^2-81n-9na^3+27a^3n^2
\\&\q-81n^3+27a^3n^3 -108n^4a^4-81an^5+108n^4a^2\mod{3^{{\rm
ord}_3n+5}}.
\endalign$$ To see the result, we note that
$$\align &a^2\e 1\mod 3,\ a^4=(a^2-1)^2+2a^2-1\e 2a^2-1\mod 9,\\&
a^5=a(a^2-1)^2+2a^3-a\e 2a^3-a\mod 9\endalign$$ and
$$\align &-81an^5+(-108a^4+108a^2)n^4+(-108a^5+27a^3)n^3
\\&\e -81an^3+27(1-a^2)n^2+(54a^3
+108a)n^3\mod {3^{{\rm ord}_3n+5}}.\endalign$$
\section*{5. A congruence for $E_{2n,a}\pmod{5^{{\rm ord}_5n+4}}$}
\pro{Theorem 5.1} Let $a$ be a nonzero integer and $n\in\Bbb N$ with
$n\ge 2$. Then
 $$\align E_{2n,a}&\e\f{2(1+(-1)^n)a^2-4a}
 {5a^2-10a+4}\\&\q+(-125a^4+250a^3(-1)^n
+250a^2(-1)^n+250a((-1)^n-1))n^3
\\&\q+((150(-1)^n+100)a^7+300a^6+(275(-1)^n
-25)a^4-25a^3(-1)^n\\&\q-25a^2(-1)^n-125a(1+(-1)^n))n^2+((-200(-1)^n-300)a^7
\\&\q+(25-200(-1)^n)a^6-(275(-1)^n+100)a^5+(30(-1)^n+105)a^4\\&\q+270a^3(-1)^n
-290a^2(-1)^n+250a((-1)^n-1))n\pmod{5^{{\rm ord}_5n+4}}.\endalign$$
\endpro
Proof. By Theorem 2.6,
$$\align E_{2n,a}&=\f1{(5(a-1)^2-1)^2}\Big\{2\cdot4^{2n}a^4+(3a^4-8a^3+4a^2)2^{2n+1}
\\&\q+4(a^4-a^3)3^{2n}+8(a^4-5a^3+6a^2-2a)-a^5\s_{k=1}^n\b{2n}{2k}5^{2k}E_{2n-2k}
\Big\}.\tag 5.1\endalign$$ As $$\align
&4^{2n}=(1+15)^n=\sum_{k=0}^n\b nk15^k\e 1+15n\mod{25},
\\&2^{2n}=(5-1)^n=\sum_{k=0}^n\b nk5^k(-1)^{n-k}\e
(-1)^n(1-5n)\mod{25},
\\&3^{2n}=(10-1)^n=\sum_{k=0}^n\b nk10^k(-1)^{n-k}\e
(-1)^n(1-10n)\mod{25},
\\&\f 1{(5a^2-10a+4)^2}=\f 1{(5(a-1)^2-1)^2}=
\f {(5(a-1)^2+1)^2}{(25(a-1)^4-1)^2}\e 10(a-1)^2+1\mod{25},
\endalign$$
we see that
$$\align
E_{2n,a}&\e(10(a-1)^2+1)\big(2a^4(1+15n)+2(3a^4-8a^3+4a^2)(-1)^n(1-5n)
\\&\q+4(a^4-a^3)(-1)^n(1-10n)+8(a^2-a)(a^2-4a+2)\big)
\\&\e2a^4(10(a-1)^2+1)+30a^4n+(6a^4-16a^3+8a^2)(-1)^n(10(a-1)^2+1-5n)
\\&\q+4(a^4-a^3)(-1)^n(10(a-1)^2+1-10n)
\\&\q+8(a^2-a)(a^2-4a+2)(10(a-1)^2+1)
\\&\e5a^2((1+(-1)^n)a^2-(-1)^n(a-2))n-10(1+(-1)^n)a^4-5(2+(-1)^n)a^3
\\&\q-2(1+6(-1)^n)a^2-a\mod{25}\endalign$$
and so
$$\align E_{2n-2,a}&\e 5a^2((1-(-1)^n)a^2+(-1)^n(a-2))(n-1)
\\&\q-10(1-(-1)^n)a^4-5(2-(-1)^n)a^3
-2(1-6(-1)^n)a^2-a\mod{25}.\endalign$$ This together with (5.1)
yields
$$\aligned E_{2n,a}&\e \f1{(5(a-1)^2-1)^2}\Big\{2\cdot4^{2n}a^4+(3a^4-8a^3+4a^2)2^{2n+1}
\\&\q+4(a^4-a^3)3^{2n}+8(a^4-5a^3+6a^2-2a)
\\&\q-a^5\b{2n}25^2(5a^2((1-(-1)^n)a^2+(-1)^n(a-2))(n-1)
\\&\q-10(1-(-1)^n)a^4 -5(2-(-1)^n)a^3
 -2(1-6(-1)^n)a^2-a)\Big\}\pmod{5^{{\rm ord}_5n+4}}.\endaligned\tag 5.2$$
 It is clear that
 $$\align &\f 1{(5(a-1)^2-1)^2}\\&=(1-5(a-1)^2)^{-2}
 =\sum_{k=0}^{\infty}\b{-2}k(-5(a-1)^2)^k
 \e \sum_{k=0}^3\b{-2}k(-5(a-1)^2)^k\\&
 \e 1+10(a-1)^2+75(a-1)^4-125(a-1)^6\pmod{5^4}.\endalign$$
 Also,
 $$\align4^{2n}&=\sum_{k=0}^n\b nk 15^k=1+n\sum_{k=1}^n\b {n-1}{k-1} \f{15^k}k
 \\&\e
 1+15n+225\f{n(n-1)}2+15^3\f{n(n-1)(n-2)}6\pmod{5^{{\rm ord}_5n+4}}\endalign$$
 $$\align2^{2n}&=\sum_{k=0}^n\b nk 5^k(-1)^{n-k}=(-1)^n+n\sum_{k=1}^n\b {n-1}{k-1}
 \f{5^k}k(-1)^{n-k}
 \\&\e
 (-1)^n\Big(1-5n+25\f{n(n-1)}2-125\f{n(n-1)(n-2)}6\Big)\pmod{5^{{\rm ord}_5n+4}}\endalign$$
 $$\align3^{2n}&=\sum_{k=0}^n\b nk 10^k(-1)^{n-k}=(-1)^n+n\sum_{k=1}^n\b {n-1}{k-1}
 \f{10^k}k(-1)^{n-k}
 \\&\e
 (-1)^n\Big(1-10n+50n(n-1)-500\f{n(n-1)(n-2)}3\Big)\pmod{5^{{\rm ord}_5n+4}}\endalign$$
 Now, from (5.2) and
the above we deduce the result.

\end{CJK*}
\end{document}